\documentclass[10pt]{article}
\usepackage{mathrsfs}
\usepackage{extarrows}
\usepackage{amsmath,amscd,amsthm,amsfonts}
\usepackage{graphicx}
\oddsidemargin=0.5cm\topmargin=-1cm \numberwithin{equation}{section}

\theoremstyle{plain}
\newtheorem{theorem}{theorem}[section]
\newtheorem{corollary}[theorem]{Corollary}

\newtheorem{fact}[theorem]{Fact}

\newtheorem{remark}[theorem]{Remark}

\newtheorem*{ThmA}{Theorem A}
\newtheorem*{ThmB}{Theorem B}

\def\disp{\displaystyle}

\begin{document}
\vskip 5.1cm
\title{Quasi-Shadowing and Quasi-Stability for Dynamically Coherent Partially Hyperbolic Diffeomorphisms\footnotetext{\\
\emph{ 2000 Mathematics Subject Classification}:37D30, 37C15, 37C50\\
\emph{Keywords and phrases}: partial hyperbolicity;
 quasi-shadowing; quasi-stability; dynamically coherent; plaque expansive.\\
 The second author is supported by NSFC(No:11001284) and Natural Science Foundation Project of CQCSTC (No:cstcjjA00003).
The third author is the corresponding author, and is supported by NSFC(No:11371120), NCET(No:11-0935), NSFHB(No:A2014205154) and the Plan of Prominent Personnel Selection and Training for the Higher Education Disciplines in Hebei Province (No:BR2-219).}}
\author {Huyi Hu$^1$, Yunhua Zhou$^2$ and Yujun Zhu$^3$   \\
\small {1. Department of Mathematics,}\\
\small {Michigan State University, East Lansing, MI 48824,USA}\\
\small {2. College of Mathematics and Statistics,}\\
\small{ Chongqing University, Chongqing, 401331.,  P. R. China}\\
\small {3. College of Mathematics and Information Science,}\\
\small {Hebei Normal University, Shijiazhuang, 050024, P.R.China}}
\date{}
\maketitle
\begin{center}
\begin{minipage}{130mm}
{\bf Abstract}:
Let $f$ be a partially hyperbolic diffeomorphism. $f$ is called has the quasi-shadowing property if for any pseudo orbit $\{x_k\}_{k\in \mathbb{Z}}$, there is a sequence  $\{y_k\}_{k\in \mathbb{Z}}$ tracing it in which $y_{k+1}$ lies in the local center leaf of $f(y_k)$ for any $k\in \mathbb{Z}$. $f$ is called topologically quasi-stable if for any homeomorphism $g$ $C^0$-close to $f$, there exist a continuous map $\pi$ and a motion $\tau$ along the center foliation such that $\pi\circ g=\tau\circ f\circ\pi$.
In this paper we prove that if $f$ is dynamically coherent then it has quasi-shadowing and topological quasi-stability properties.

\end{minipage}
\end{center}

\section{Introduction}

The goal of this paper is to investigate the quasi-shadowing property and its relation to the quasi-stability for
partially hyperbolic diffeomorphisms under certain assumptions on the foliations.

It is well known that any Anosov diffeomorphism $f$ on a closed manifold $M$ has the shadowing property (see \cite{Anosov} and \cite{Bowen} for example). That is, for any $\varepsilon>0$ there exists $\delta>0$ such
that every pseudo orbit  $\xi=\{x_{k}\}_{-\infty}^{+\infty}$ with $
\sup_{k\in \mathbb{Z}} d(f(x_{k}),\;x_{k+1})\leq \delta$ is $\varepsilon$-shadowed ($\varepsilon$-traced) by some
true orbit orb$(x)$ in the sense of $
\sup_{k\in \mathbb{Z}} d(f^k(x),\;x_{k})\leq \varepsilon$. However,  for any
 partially hyperbolic diffeomorphism  we cannot expect that in general
the shadowing property holds since in this case a center direction is allowed
in addition to the hyperbolic directions. Therefore, how to find an analogous property is interesting.

Recently,  we showed that any partially hyperbolic diffeomorphism $ f:M\longrightarrow M$  without any additional assumption has the \emph{quasi-shadowing property} in the following sense: for any pseudo orbit $\{x_k\}_{k\in \mathbb{Z}}$, there is a sequence  $\{y_k\}_{k\in \mathbb{Z}}$ tracing it
in which $y_{k+1}$ is obtained from $f(y_k)$ by a motion $\tau$
along center direction.  In particular, if $f$ has $C^1$ center foliation, the motion $\tau$ can be chosen along the center foliation and hence has more geometrical meaning(\cite{HZZ}). To obtain the above quasi-shadowing property, we used a  unified ``analytic" method which is a modification of that in \cite{Anosov} and \cite{Kato}, where it was used to get the shadowing and stability properties for Anosov diffeomorphisms and Anosov flows respectively. In this paper, we shall obtain the quasi-shadowing property with the motion $\tau$ also moving points along the center foliation under the dynamical coherence assumption which is weaker than smoothness. In this case, the  above ``analytic" method seems ineffective, we turn to  a ``geometric" method which relies on more detailed information about the foliations of the system. This result is
given in Theorem A. We mentioned that Kryzhevich and Tikhomirov \cite{Tikhomirov}  got a similar result recently using a different strategy which combines the  ``analytic" and ``geometric" techniques. Actually, in our proof we borrow the elegant idea of Bowen (\cite{Bowen}) to construct the quasi-shadowing sequence for a pseudo orbit and hence in some sense this construction is more intuitive.

Shadowing property is a powerful tool to investigate the stability property of a system. In  \cite{Walters1},  Walters showed that shadowing property implies topological stability for any expansive homeomorphism on a compact metric space. Hence any Anosov diffeomorphism $f:M\longrightarrow  M$ is  topologically stable (\cite{Walters}), that is, for any homeomorphism $g$ $C^0$-close to $f$ there exists a continuous map $\pi$ onto $M$ such that
\begin{equation}\label{conjugacy}
\pi\circ g=f\circ\pi.
\end{equation}
For partially hyperbolic diffeomorphisms,
we can not expect such stabilities because of the existence
of the center direction.

Recently,  Hu and Zhu \cite{HZ} showed that any partially hyperbolic diffeomorphism $ f:M\longrightarrow M$  without any additional assumption is \emph{topologically quasi-stable} in the following sense: for any homeomorphism $g$ $C^0$-close to $f$, there exist a
continuous map $\pi$ from $M$ to itself and a family of locally
defined maps $\{\tau_x:x\in M\}$, which map points
along the center direction, such that
\begin{equation}\label{conjugacy1}
\pi\circ g(x)=\tau_{f(x)}\circ f\circ\pi(x) \quad  \text{for all} \
x\in M.
\end{equation}  In particular, if $f$ has $C^1$ center foliation, the motion $\tau$ can be chosen along the center foliation. To get these quasi-stability results, we also used a unified ``analytic" method. However, if one want to get the quasi-stability with the motion $\tau$ also moving points along the center foliation without the smoothness  assumption on the center foliation, this ``analytic" method also seems ineffective. Fortunately, now we can resort to the quasi-shadowing we obtain in Theorem A.  Recently, Kryzhevich \cite{Kryzhevich} showed that dynamical coherence implies plaque expansivity. Quasi-shadowing and plaque expansivity for partially hyperbolic diffeomorphisms are the counterparts of shadowing and expansivity respectively for Anosov diffeomorphisms. Hence,
we can apply these two properties to get the topological quasi-stability with the desired motion $\tau$ for the partially hyperbolic diffeomorphism under the assumption of dynamical coherence (Theorem B).

\section{Definition and statement of results}\label{DSN}

Everywhere in this paper, we assume that $M$ is a smooth
$m$-dimensional compact Riemannian manifold. We denote by $\|\cdot\|$
and $d(\cdot,\cdot)$ the norm on $TM$ and the metric
on $M$ induced by the Riemannian metric, respectively.

A diffeomorphism $f: M\to M$ is said to be
\emph{(uniformly) partially hyperbolic} if there exist
numbers $\lambda,\lambda',\mu$ and $\mu'$ with
$0<\lambda<1<\mu$ and $\lambda<\lambda'\leq \mu'<\mu$,
and an invariant decomposition $T_xM=E_x^s\oplus E_x^c\oplus E_x^u$
\ $\forall x\in M$, such that for any $n\ge 0$,
\begin{eqnarray*}
\|d_xf^nv\|  &\!\!\!\!\!  \le C\lambda^n\|v\|\ \   &\text{as} \ v\in E^s(x), \\
C^{-1}(\lambda')^n\|v\| \le \|d_xf^nv\|  &\!\!\!  \le C(\mu')^n\|v\|
               \ &\text{as} \   v\in E^c(x),\\
C^{-1}\mu^n\,\|v\| \leq  \|d_xf^nv\|&\ &\text{as} \  v\in E^u(x)
\end{eqnarray*}
hold for some number $C>0$.  The subspaces
$E_x^s, E_x^c$ and $E_x^u$ are called
\emph{stable, center} and \emph{unstable} subspace, respectively.
Via a change of Riemannian metric we always assume that $C=1$.
Moreover, for simplicity of the notation, we assume that
$\disp \lambda=\frac{1}{\mu}$. For general theory of partially hyperbolic system, we refer to \cite{Hirsch},
\cite{Pesin}, \cite{Barreira} and \cite{Bonatti}.

A partially hyperbolic diffeomorphism $f$ is called \emph{dynamically coherent} if
$E^{cu}:=E^c\oplus E^u$, $E^{c}$, and $E^{cs}:=E^c\oplus E^s$ are
integrable, and everywhere tangent to $\mathcal{W}^{cu}$,
$\mathcal{W}^{c}$ and $\mathcal{W}^{cs}$, called the {\it
center-unstable}, \emph{center} and \emph{center-stable foliations},
respectively; and $\mathcal{W}^c$ and $\mathcal{W}^u$ are
subfoliations of $\mathcal{W}^{cu}$, while $\mathcal{W}^c$ and
$\mathcal{W}^s$ are subfoliations of $\mathcal{W}^{cs}$. For any $x\in M$, denote $W^\omega(x), \omega\in \{s,c,u,cs,cu\}$, the stable, center, unstable, center-stable, center-unstable manifold of $x$ respectively. For $\varepsilon>0$, denote $W^\omega_{\varepsilon}(x), \omega\in \{s,c,u,cs,cu\}$ the corresponding local manifolds of size $\varepsilon$. We also call these local manifolds \emph{plaques}.

\begin{ThmA}\label{theoremA}
Let $f:M \to M$ be a partially hyperbolic diffeomorphism.
If $f$ is dynamically coherent then $f$ has the
quasi-shadowing property in the following sense:  for any $\varepsilon>0$ there exists $\delta>0$ such that for any $\delta$-pseudo
orbit $\{x_k\}_{k\in \mathbb{Z}}$ of $f$, there exists a sequence
$\{y_k\}_{k\in \mathbb{Z}}$ with $y_k\in W^c_{\varepsilon}(f(y_{k-1}))$
such that $d(x_k,y_k)<\varepsilon$.

\end{ThmA}

Now we consider an application of our result to the quasi-stability.  $f$ is called \emph{plaque expansive} with respect to the center foliation  $\mathcal{W}^c$ if there exists $\eta>0$ such that for any $\eta$-pseudo orbits $\{x_n\}_{n=-\infty}^{\infty}$ and $\{y_n\}_{n=-\infty}^{\infty}$ in which $f(x_n)$ and $f(y_n)$  lie in the center plaque $W^c_{\eta}(x_{n+1})$ and $W^c_{\eta}(y_{n+1})$ respectively, $x_n$ and $y_n$ must lie in a common center plaque, i.e., $x_n\in W^c_{\eta}(y_n)$. Such an $\eta$ is called a \emph{plaque expansiveness constant} of $f$. Recently, Kryzhevich \cite{Kryzhevich} showed that for any partially hyperbolic diffeomorphism,  dynamical coherence implies plaque expansivity. Applying quasi-shadowing and plaque expansivity properties, we can get the following quasi-stability property.

\begin{ThmB}\label{ThmB}
Let $f:M \to M$ be a partially hyperbolic diffeomorphism.
If $f$ is dynamically coherent, then $f$ has the topological quasi-stability in the following sense:  for any
$\varepsilon>0$ there exists $\delta>0$ such that
for any homeomorphism $g$ of $M$ with $d(f,g)<\delta$
there exist a  surjective continuous map $\pi:M\to M$ and a family of motions $\{\tau_x:x\in M\}$  which move points along the center foliation and are continuously dependent on $x$ such that
\begin{equation}\label{conjugacy2}
\pi\circ g(x)=\tau_x\circ f\circ \pi(x),\;\;x\in M,
\end{equation}
i.e., the following commutative diagram
$$
\begin{CD}
{\quad}& \quad &M&  &\xlongrightarrow[]{\quad\quad\quad\quad\quad g\quad\quad\quad\quad}&  &M&  \\
 \quad& \pi&@VVV&     &@VVV&\pi &\quad&\\
{\quad}&\quad &M&
&\xlongrightarrow{\quad\quad f\quad\quad}  M
\xlongrightarrow{\tau_{(\cdot)}}& &M&
\end{CD}\qquad
$$
holds.
\end{ThmB}

\begin{remark}
In Theorem B, the family of motions $\{\tau_x:x\in M\}$ are continuously dependent on $x$ means that  $\tau_x\circ f\circ \pi(x)$ tends to $\tau_y\circ f\circ \pi(y)$ as $y$ tends to $x$.
Here we write the motions in the form $\tau_x, x\in M$ to indicate their dependence on $x$. In fact, from the definitions of $\pi$ and $\tau_x$ in the proof of this theorem, there possibly exist two points $x,y\in M$ with $\pi(x)=\pi(y)$ but  $\pi g(x)\neq\pi g(y)$, and hence $\tau_x\circ f\circ \pi(x)\neq\tau_y\circ f\circ \pi(y)$. However, when the homeomorphism $g$ in Theorem B is also expansive, then the family of motions  $\{\tau_x:x\in M\}$ is ``coherent" and form a homeomorphism.
\end{remark}

\begin{corollary}
Let $f$ be as in Theorem B.  Then for any expansive homeomorphism $g$ $C^0$ sufficiently close to $f$
there exist two homeomorphisms $\pi$ and $\tau$ on $M$ such that
\begin{equation}\label{feqncoro}
\pi\circ g=\tau\circ f\circ \pi,
\end{equation}
where $\tau$ maps every point to its image along the center foliation.
\end{corollary}

\section{Quasi-shadowing}

 To prove Theorem A, we first list some facts about the foliations. By the dynamical coherence and transversality of the foliations, we have the following fact.

 \begin{fact}
 There exists $\delta_0>0,L_0>1$ such that for any $0<\delta<\delta_0$ and $x,y\in M$ satisfying $d(x,y)<\delta$, the following statements hold.

 1) The intersection $W^{cu}_{L_0\delta}(x)\bigcap W^{s}_{L_0\delta}(y)$ consists of a single point and so does the intersection $W^{cs}_{L_0\delta}(x)\bigcap W^{u}_{L_0\delta}(y)$.

 2) When $x, y$ lie in a common center-unstable (resp. center-stable) leaf, the intersection $W^{c}_{L_0\delta}(x)\bigcap W^{u}_{L_0\delta}(y)$ (resp. $W^{c}_{L_0\delta}(x)\bigcap W^{s}_{L_0\delta}(y)$) consists of a single point.
  \end{fact}

 Fix $x\in M$, $0<r<\frac{L_0\delta_0}{3}$ and consider the family of center plaques
$$
\mathcal{L}^{c,u}_r(x)=\{W^c_r(z): z\in B^{cu}(x, r)\},
$$
where $B^{cu}(x, r)=B(x, r)\cap W^{cu}(x)$. Clearly, the union of the elements in
$\mathcal{L}^{c,u}_r(x)$ form a center-unstable plaque, for convenience, we also denote it by $\mathcal{L}^{c,u}_r(x)$. Choose two unstable plaques $D^1$ and $D^2$ which are transverse to $\mathcal{L}^{c,u}_r(x)$ and define a map
$h: D^1\rightarrow D^2$ by setting
$$
h(y)=D^2\cap W^c_r(y)\;\;\; \text{for}\;\;\;y\in D^1.
$$
Clearly, the map is well defined and it is a
homeomorphism onto its image. It is called an unstable \emph{holonomy map} along the center leaves inside
a center-unstable plaque. Similarly, we can define a family of center plaques
$\mathcal{L}^{c,s}_r(x)$ inside a local center-stable manifold of $x$ and define a stable
holonomy map along the center leaves for any two stable plaques $D^1$ and $D^2$ which are transverse to $\mathcal{L}^{c,s}_r(x)$. Since $M$ is
compact, these holonomy maps have the locally equivalent continuity in
the following sense.

 \begin{fact}
 For any $\alpha>0$, there exist $0<r_1,
r_2<\frac{L_0\delta_0}{3}$ such that for any $x\in M$ and any two unstable (resp. stable) plaques $D^1$
and $D^2$ which are transverse to $\mathcal{L}^{c,u}_{r_1}(x)$ (resp. $\mathcal{L}^{c,s}_{r_1}(x)$), the corresponding unstable (resp. stable) holonomy map $h:
D^1\rightarrow D^2$ along the center leaves has the following property: for any $z, z'\in
D^1$
\begin{equation}\label{holonomy}
d(z, z')<r_2\;\;\; \text{implies}\;\;\;d(h(z), h(z'))<\alpha.
\end{equation}

\end{fact}

Note that for any $k\in \mathbb{Z}$ the corresponding foliations $\mathcal{W}^{\omega}, \omega\in \{s, c, u, cs, cu\}$ of $f^k$ coincide with that of $f$, then we have the following fact.
 \begin{fact}
$f$ has the quasi-shadowing property if and only if $f^k$  so does.
\end{fact}
This fact ensure me to assume the hyperbolicity constant $\lambda$ is sufficiently small to meet our needs.

Take $\delta_1>0$ such that
\begin{equation}\label{constant}
\begin{split}
&d(f(x),f(y))<\lambda d(x,y) \;\;\;\mbox{for} \;\;\;y\in W^s_{\delta_1} (x),\\
&d(f^{-1}(x),f^{-1}(y))<\lambda d(x,y) \;\;\;\mbox{for} \;\;\;y\in W^u_{\delta_1} (x).
\end{split}
\end{equation}

\begin{proof}[Proof of Theorem A]

For any $\varepsilon>0$, take $\alpha<\frac{\varepsilon}{3}$
(correspondingly, $r_1$ and $r_2$ in Fact 2.2 are taken with respect
to $\alpha$), and then take $\delta$ and $\lambda$ all small enough
satisfying
\begin{equation}\label{cond0}
2\lambda L_0<1,
\end{equation}
\begin{equation}\label{cond1}
\delta(1+2L_0+2\lambda L_0)<\frac{\varepsilon}{3}
\end{equation}
and
\begin{equation}\label{cond2}
\lambda(2L_0\delta+\alpha)<r_2,
\end{equation}
and such that the following two statements hold.

\textbf{Statement 1}.  Any stable, center, unstable, center-stable and center-unstable plaques which will be concerned in the following proof are all with the size less than $\max\{L_0\delta_0, \delta_1\}$, hence 1) and 2) in Fact 2.1 and the inequalities in  (\ref{constant}) hold.

\textbf{Statement 2}. Any center-stable and center-unstable plaques which will be concerned in the following proof are contained in some  $\mathcal{L}^{c,s}_{r_1}$ and $\mathcal{L}^{c,u}_{r_1}$  respectively, hence Fact 2.2 holds.

We emphasize that the conditions (\ref{cond0}), (\ref{cond1}), (\ref{cond2}) are a little bit technical. In fact, you can first ignore them and the estimations of the sizes of many plaques and the distance between the points we have to deal with in the following and quickly catch the idea of how the quasi-shadowing sequences are constructed.

Let $\xi=\{x_i\}_{i=-\infty}^{+\infty}$ be a $\delta$-pseudo orbit of $f$. We will find a sequence $\{y^*_i\}_{i=-\infty}^{+\infty}$ which $\varepsilon$-quasi-shadows $\xi$ in three steps.

\emph{Step 1}.  Find a sequence $\{y^u_i\}_{i=0}^{+\infty}$ which
$\frac{2\varepsilon}{3}$-quasi-shadows the positive half sequence
$\{x_i\}_{i=0}^{+\infty}$.

Firstly, we consider a finite piece $\{x_i\}_{i=0}^n$ of $\xi$ for
$n\geq 1$. In the following, we will define four sequences
$\{z_i\}_{i=1}^n$, $\{z_i'\}_{i=1}^{n}$ ,  $\{y_i\}_{i=0}^{n-1}$ and
$\{y_i'\}_{i=1}^{n-1}$ successively. The existence and uniqueness of
these points are ensured by Statement 1 and Statement 2. The first
two sequences $\{z_i\}_{i=1}^n$ and $\{z_i'\}_{i=1}^{n}$
are defined as follows. Since $d(f(x_0), x_1)<\delta$, from Fact 2.1, take
$$
\{z_1\}=W_{L_0\delta}^{cu}(f(x_0))\cap W_{L_0\delta}^{s}(x_1)\;\;\;\text{and}\;\;\;
\{z_1'\}=W_{L_0\delta}^{u}(f(x_0))\cap W_{L_0\delta}^{cs}(x_1).
$$
Note that by (\ref{constant}) and (\ref{cond0}),
$$
d(f(z_1), x_2)\leq  d(f(z_1), f(x_1))+d(f(x_1), x_2)
< \lambda L_0\delta+\delta< 2\delta.
$$
Let
$$\{z_2\}=W_{2L_0\delta}^{cu}(f(z_1))\cap
W_{2L_0\delta}^{s}(x_2)\;\;\;\text{and}\;\;\;\{z_2'\}=W_{2L_0\delta}^{u}(f(z_1))\cap W_{2L_0\delta}^{cs}(x_2).
$$
Now assume that for any $2\leq k\leq n-1,$
$$
\{z_k\}=W_{2L_0\delta}^{cu}(f(z_{k-1}))\cap W_{2L_0\delta}^{s}(x_k)\;\;\;\text{and}\;\;\;\{z_k'\}=W_{2L_0\delta}^{u}(f(z_{k-1}))\cap
W_{2L_0\delta}^{cs}(x_k).
$$
 Then
$$
d(f(z_k), x_{k+1})\leq  d(f(z_k), f(x_k))+d(f(x_k),
x_{k+1})<2\lambda L_0\delta+\delta <2\delta.
$$
Hence we can take
$$\{z_{k+1}\}=W_{2L_0\delta}^{cu}(f(z_{k}))\cap
W_{2L_0\delta}^{s}(x_{k+1})\;\;\;\text{and}\;\;\;\{z_{k+1}'\}=W_{2L_0\delta}^{u}(f(z_{k}))\cap
W_{2L_0\delta}^{cs}(x_{k+1}).
$$
By the dynamical coherence, for any
$1\leq i\leq n$, $z_i$ and $z_i'$ lie in a common center plaque. Now
we define the last two sequences $\{y_i\}_{i=0}^{n-1}$ and
$\{y_i'\}_{i=1}^{n-1}$ as follows. Let $y_{n-1}=f^{-1}(z_n')$ and
$y_{n-1}'$ be the unique point in
the intersection of the center plaque of $y_{n-1}$ and the unstable
plaque of $z_{n-1}'$. Assume for any $2\leq k\leq n-1$, $y_k$ and $y_k'$ are
defined, then take $y_{k-1}=f^{-1}(y_k')$ and $y_{k-1}'$ to be the
unique point in the intersection of the center plaque of $y_{k-1}$ and
the unstable plaque of $z_{k-1}'$.  Finally let $y_0=f^{-1}(y_1')$.

From the definition of these four sequences and statement 2, we can
see that each element in the family
$$
\big\{\{x_0, y_0\}, \{f(x_0),
z_1', y_1'\}, \{f(z_1), z_2', y_2'\}, \cdots,  \{f(z_{n-2}),
z_{n-1}', y_{n-1}'\}, \{f(z_{n-1}), z_{n}'\}\big\}
$$
lies in an
unstable plaque, and each element in the family $\big\{\{y_i, z_i,
y_i', z_i'\}\big\}_{i=1}^{n-1}$ lies in a center-unstable plaque
$\mathcal{L}_{r_1}^{c,u}$. Moreover, $\{y_i\}_{i=0}^{n-1}$
quasi-shadows $\{x_i\}_{i=0}^{n-1}$. Now we estimate, for $1\leq i\leq n-1$,  the distance
between $x_i$ and $y_i$. Since $z_{n}'\in
W_{2L_0\delta}^{u}(f(z_{n-1}))$ and $y_{n-1}=f^{-1}(z_n')$, by
(2.2), we have
$$
d(z_{n-1}, y_{n-1})<\lambda d(f(z_{n-1}), z_{n}')\leq
2\lambda L_0\delta<r_2.  \;\;\;(\mbox{by}\;(\ref{cond2}))
$$
So
\begin{eqnarray}
d(x_{n-1}, y_{n-1})&\leq&d(x_{n-1}, z_{n-1})+d(z_{n-1}, y_{n-1})\notag\\
&<&2 L_0\delta+2\lambda
L_0\delta<\frac{\varepsilon}{3},\;\;\;(\mbox{by}\;(\ref{cond1}))\notag
\end{eqnarray}
and $d(z_{n-1}',y_{n-1}')<\alpha$ by Fact 2.2 . Then
$$
d(f(z_{n-2}), y_{n-1}')\leq d(f(z_{n-2}), z_{n-1}')+d(z_{n-1}',
y_{n-1}')<2L_0\delta+\alpha.
$$
Hence,
$$
d(z_{n-2}, y_{n-2})<\lambda
d(f(z_{n-2}), y_{n-1}')<\lambda(2L_0\delta+\alpha), \;\;\;(\mbox{by}\;(\ref{constant}))
$$
and
\begin{eqnarray}
d(x_{n-2}, y_{n-2})&\leq& d(x_{n-2}, z_{n-2})+d(z_{n-2},
y_{n-2})\notag\\&<&2L_0\delta+\lambda(2L_0\delta+\alpha)<\frac{2\varepsilon}{3}.
\;\;\;(\mbox{by}\;(\ref{cond1}))\notag
\end{eqnarray}
Inductively, we can get that for $3\leq i\leq n-1$,
$$
d(x_{n-i},y_{n-i})\leq d(x_{n-i},z_{n-i})+d(z_{n-i},y_{n-i})< 2L_0\delta+\lambda(2L_0\delta+\alpha)<\frac{2\varepsilon}{3}
$$
and
$$
d(x_0, y_0)<\lambda d(f(x_0),
y_1')\leq \lambda [d(f(x_0),z_1')+d(z_1',y_1')]<\lambda(L_0\delta+\alpha)<\frac{2\varepsilon}{3}.
$$
Therefore, we prove that $\{y_i\}_{i=0}^{n-1}$
$\frac{2\varepsilon}{3}$-quasi-shadows $\{x_i\}_{i=0}^{n-1}$. In
fact, by the construction, we can see that the sequence
$\{y_i\}_{i=0}^{n-1}$ is uniquely determined by $\{x_i\}_{i=0}^n$,
and $y_0\in W^u(x_0), y_i\in W^{cu}(f^i(x_0)) $ for any $1\leq i\leq
n-1$. For convenience, we relabel $\{y_i\}_{i=0}^{n-1}$ by
$\{y_{i,n}\}_{i=0}^{n-1}$ to indicate its dependence on $n$. Let
$y^u_0$ be one limit point of $\{y_{0,n}\}_{n=0}^\infty$. Obviously,
$y^u_0\in W^u(x_0)$. Now we define two sequences $\{(y_{i}^{u})'\}_{i=0}^\infty$ and
$\{y_{i}^{u}\}_{i=0}^\infty$ successively as follows. Let $(y_1^u)'=f(y^u_0)$ and $y^u_1$ be
the unique point in the intersection of the unstable plaque of $z_1$
and the center plaque of $(y_1^{u})'$. Assume that for any $1\leq
i\leq k, (y_i^u)'$ and $y^u_i$ are defined, then we take $(y_{k+1}^u)'=f(y^u_k)$ and $y^u_{k+1}$ to be the
unique point in the intersection of the unstable plaque of $z_{k+1}$
and the center plaque of $(y_{k+1}^u)'$. It is easy to see that
$y^u_i\in W^{cu}(f^i(x_0)) $ for any $i\geq 1$ and the sequence
$\{y_{i}^{u}\}_{i=0}^\infty$ $\frac{2\varepsilon}{3}$-quasi-shadows
$\{x_i\}_{i=0}^\infty$.

\emph{Step 2}.  Find a sequence $\{y^s_i\}_{i=-\infty }^{0}$ which
$\frac{2\varepsilon}{3}$-quasi-shadows the negative half sequence
$\{x_i\}_{i=-\infty }^{0}$.

Since the strategy in this step is similar to that in step 1 except for the type
of plaques in which the quasi-shadowing sequence lies, we only give
the outline of the construction of $\{y_i^s\}_{i=-\infty }^{0}$. For
simplicity, we assume that for any $i\leq 0$,
$d(f^{-1}(x_i),x_{i-1})<\delta$ (otherwise, we can take
$0<\delta'<\delta$ such that for any $\delta'$-pseudo orbit
$\{x_i\}_{i=-\infty}^{+\infty}$, we have
$d(f^{-1}(x_i),x_{i-1})<\delta$ for  $i\leq 0$ and show that
$\{x_i\}_{i=-\infty}^{+\infty}$ can be $\varepsilon$-quasi-shadowed
by some sequence of points). Firstly, for any finite piece
$\{x_i\}_{i=n}^0 (n\leq -1)$ of $\xi$, we define four sequences $\{z_i\}_{i=n}^{-1}, \{z_i'\}_{i=n}^{-1}$,
$\{y_i\}_{i=n+1}^0$ and $\{y_i'\}_{i=n+1}^{-1}$ successively as
follows. Let
$$
z_{-1}= W^{s}_{L_0\delta}(f^{-1}(x_0))\cap
W^{cu}_{L_0\delta}(x_{-1}),\;\;z_{-1}'=
W^{cs}_{L_0\delta}(f^{-1}(x_{0}))\cap W^{u}_{L_0\delta}(x_{-1});
$$
$$
z_{-2}= W^{s}_{2L_0\delta}(f^{-1}(z_{-1}'))\cap
W^{cu}_{2L_0\delta}(x_{-2}),\;\;z_{-2}'=
W^{cs}_{2L_0\delta}(f^{-1}(z_{-1}'))\cap W^{u}_{2L_0\delta}(x_{-2});
$$
$$
\cdots\cdots
$$
$$
z_{n}= W^{s}_{2L_0\delta}(f^{-1}(z_{n+1}'))\cap
W^{cu}_{2L_0\delta}(x_{n}),\;\;z_{n}'=
W^{cs}_{2L_0\delta}(f^{-1}(z_{n+1}'))\cap
W^{u}_{2L_0\delta}(x_{n}).
$$
Let $y_{n+1}'=f(z_n)$, $y_{n+1}$ be the unique point in the
intersection of the center plaque of $y_{n+1}'$ and the stable
plaque of $z_{n+1}$. Inductively define $y_i'=f(y_{i-1})$ and $y_i$ to be the unique point in the intersection
of the center plaque of $y_{i}'$ and the stable plaque of $z_{i}$ for any $n+2\leq i\leq -1$,
and let $y_0=f(y_{-1})$. From the construction and statement 2, we
can see that each element in the family
$$
\big\{\{x_0, y_0\},
\{f^{-1}(x_0), z_{-1}, y_{-1}\}, \{f^{-1}(z_{-1}'), z_{-2},
y_{-2}\}, \cdots, \{f^{-1}(z_{n+2}'), z_{n+1}, y_{n+1}\},
\{f^{-1}(z_{n+1}'), z_{n}\}\big\}
$$
lies in a stable plaque and each
element in the family $\big\{\{z_i, y_i, z_i',
y_i'\}\big\}_{i=n+1}^{-1}$
lies in a center-stable plaque
$\mathcal{L}_{r_1}^{c,s}$. Moreover $\{y_i\}_{i=n+1}^{0}$
quasi-shadows $\{x_i\}_{i=n+1}^{0}$. Now we estimate the distance
between $x_i$ and $y_i$ for any $n+1\leq i\leq 0$. Since
$f^{-1}(z_{n+1}')$ and $z_{n}$ lie in a common stable plaque and
$d(f^{-1}(z_{n+1}'), z_n)<2L_0\delta$, then
$$
d(z_{n+1}', y_{n+1}')<\lambda d(f^{-1}(z_{n+1}'), z_n)<2\lambda
L_0\delta<r_2.
$$
Hence by Fact 2.2, $d(z_{n+1}, y_{n+1})<\alpha$ and then
\begin{eqnarray}
d(x_{n+1}, y_{n+1})&\leq& d(x_{n+1}, f^{-1}(x_{n+2}))+
d(f^{-1}(x_{n+2}),
f^{-1}(z_{n+2})')\notag\\&\quad&+d(f^{-1}(z_{n+2})',
z_{n+1})+d(z_{n+1}, y_{n+1})\notag\\&<& \delta +2\lambda
L_0\delta+2L_0\delta+\alpha<\frac{\varepsilon}{3}+\frac{\varepsilon}{3}=\frac{2\varepsilon}{3}.\notag\;\;\;(\text{by}
(\ref{cond1}))
\end{eqnarray}
Similarly, we have $d(x_{i}, y_{i})<\frac{2\varepsilon}{3}$ for any
$n+2\leq i\leq -1$. Moreover,
$$
d(x_{0}, y_{0})<\lambda d(f^{-1}(x_{0}),
y_{-1})\leq\lambda [d(f^{-1}(x_{0}), z_{-1})+d(z_{-1},
y_{-1})]<\lambda
(L_0\delta+\alpha)<\frac{2\varepsilon}{3}.
$$
Therefore, we prove that $\{y_i\}_{i=n+1}^0$
$\frac{2\varepsilon}{3}$-quasi-shadows $\{x_i\}_{i=n+1}^0$. Relabel $\{y_i\}_{i=n+1}^{0}$ by $\{y_{i,n}\}_{i=n+1}^{0}$ and let
$y_0^s$ be one limit point of $\{y_{0,n}\}_{n=-\infty}^0$. Obviously
$y_0^s\in W^s(x_0)$. Now define two sequences
$\{y_{i}^s\}_{i=-\infty}^0$ and $\{(y_{i}^s)'\}_{i=-\infty}^0$ as follows. Let
$y_{-1}^s=f^{-1}(y_0^s)$, $(y_{-1}^s)'$ be the unique
point in the intersection of the stable plaque of $z_{-1}'$ and the
center plaque of $y_{-1}^s$. Inductively define $y_{i}^s=f^{-1}((y_{i+1}^s)')$  and $(y_i^s)' $ to be the  unique
point in the intersection of the stable plaque of $z_i'$ and the
center plaque of $y_i^s$ for any $i\leq -2$. Clearly,
$y_{i}^s\in W^{cs}(f^i(x_0))$ for any $i\leq -1$ and
$\{y_i^s\}_{i=-\infty}^{0}$ $\frac{2\varepsilon}{3}$-quasi-shadows
$\{x_i\}_{i=-\infty }^{0}$.

\emph{Step 3}.  Construct the desired sequence $\{y^*_i\}_{i=-\infty
}^{+\infty}$.

Note that $d(y_0^s, y_0^u)\leq d(y_0^s, x_0)+d(x_0, y_0^u)<
2\lambda (L_0\delta+\alpha)$ . So by Fact 2.1, we can take
$$
y_0^*=W_{2\lambda L_0(L_0\delta+\alpha)}^{s}(y_0^u)\cap
W_{2\lambda L_0(L_0\delta+\alpha)}^{cu}(y_0^s)\;\;\;
\text{and}\;\;\;
(y_0^*)'=W_{2\lambda L_0(L_0\delta+\alpha)}^{cs}(y_0^u)\cap
W_{2\lambda L_0(L_0\delta+\alpha)}^{u}(y_0^s).
$$
We now define two
sequences $\{y^*_i\}_{i=-\infty }^{+\infty}$ and
$\{(y^*_i)'\}_{i=-\infty }^{+\infty}$ as follows. For the positive
direction, inductively define $(y^*_i)'=f(y^*_{i-1})$ and $y^*_i$ to
be the single point in the intersection of the center plaque of
$(y^*_i)'$ and the stable plaque of $y_{i}^u$ for $i\geq 1$. It is
obvious that each element of the family $\big\{\{(y_i^u)', (y^*_i)',
y_i^u, y^*_i\}\big\}_{i=1}^{\infty}$ lies in a center-stable plaque
$\mathcal{L}_{r_1}^{c,s}$ and $d((y_i^u)', (y^*_i)')<r_2$, hence by Fact 2.2, $d(y_i^u, y^*_i)\leq
\alpha<\frac{\varepsilon}{3}$ for $i\geq 1$. Therefore, for $i\geq 1$,
\begin{eqnarray}\label{dist1}
d(x_i, y^*_i)&\leq&d(x_i, y^u_i)+d(y^u_i,
y^*_i)\leq\frac{2\varepsilon}{3}+\frac{\varepsilon}{3}=\varepsilon.
\end{eqnarray}
For the negative direction,  inductively define
$y^*_i=f^{-1}((y^*_{i+1})')$ and $(y^*_i)'$ to be the single point in
the intersection of the center plaque of $y^*_i$ and the unstable
plaque of $(y_{i}^s)'$ for $i\leq -1$. Similar to that in the positive direction, each element of
the family $\big\{\{(y_i^s)', (y^*_i)', y_i^s,
y^*_i\}\big\}_{i=-\infty}^{-1}$ lies in a center-unstable plaque
$\mathcal{L}_{r_1}^{c,u}$ and $d(y_i^s, y^*_i)\leq
\frac{\varepsilon}{3}$ for $i\leq -1$. Hence  for $i\leq -1$,
\begin{eqnarray}\label{dist2}
d(x_i, y^*_i)\leq d(x_i, y^s_i)+d(y^s_i,
y^*_i)\leq\frac{2\varepsilon}{3}+\frac{\varepsilon}{3}=\varepsilon.
\end{eqnarray}
Also note
\begin{eqnarray}\label{dist3}
d(x_0, y^*_0)\leq d(x_0, y^u_0)+d(y^u_0,
y^*_0)<\lambda(L_0\delta+\alpha)+2\lambda
L_0(L_0\delta+\alpha)<\varepsilon.
\end{eqnarray}
By (\ref{dist1}), (\ref{dist2}) and (\ref{dist3}), we conclude that the sequence
$\{y^*_i\}_{i=-\infty }^{+\infty}$ $\varepsilon$-quasi-shadows the
$\delta$-pseudo orbit $\{x_i\}_{i=-\infty}^{+\infty}$. Please see the following figure to understand the construction of the above sequences of  points.

\begin{center}
 \includegraphics[height=7.5cm]{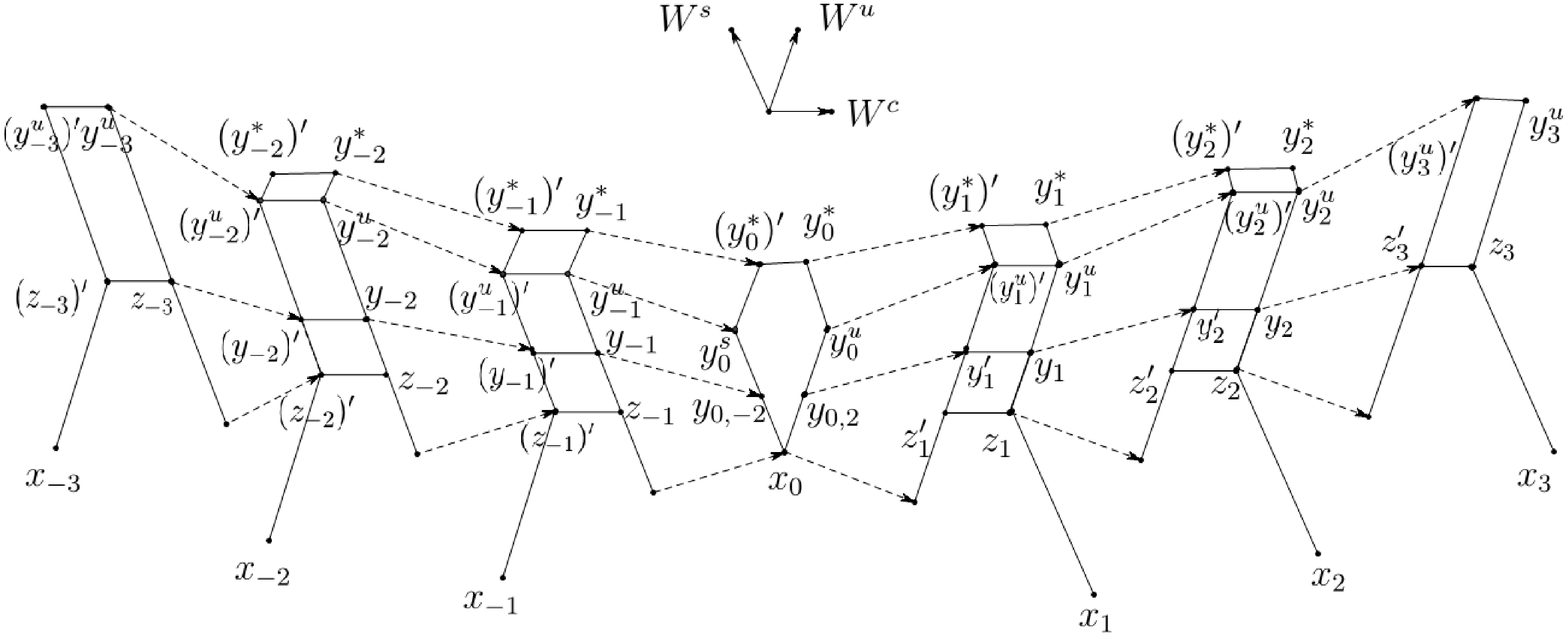}
     \label{figure}
\end{center}

This completes the proof of Theorem A.

\end{proof}

\section{Quasi-stability}

\begin{proof}[Proof of Theorem B]
Let $f:M\rightarrow M$ be a partially hyperbolic diffeomorphism
which is dynamically coherent. By \cite{Kryzhevich}, $f$ is plaque expensive with respect to the center foliation. Let $\eta$ be a plaque expensiveness constant.

Let $\varepsilon<\frac{\eta}{2}$. By Theorem A, there exists
$\delta>0$ such that for any $\delta-$pseudo orbit $\{x_i\}_{i=-\infty}^{+\infty}$ of $f$, there exists a sequence $\{y_i^*\}_{i=-\infty}^{+\infty}$ with $y_i^*\in W_{\varepsilon}^{c}(f(y_{i-1}^*))$ and
$d(x_i, y_i^*)<\varepsilon$ for any $i\in \mathbf{Z}$. Remember that in the proof of Theorem
A, $\{y_i^*\}_{i=-\infty}^{+\infty}$ is uniquely determined by
$\{y_i^s\}_{i=-\infty}^0$ and  $\{y_i^u\}_{i=0}^{+\infty}$, and
$\{y_i^s\}_{i=-\infty}^0$ and  $\{y_i^u\}_{i=0}^{+\infty}$ are
uniquely determined by $y_0^s$ and $y_0^u$ respectively. Note that
$y_0^s$ (resp. $y_0^u$) is one limit point of the fixed family
$\{y^s_{0, n}\}_{n=-\infty}^0$ (resp. $\{y^u_{0, n}\}_{n=0}^{+\infty}$)
which is uniquely determined by $\{x_n\}_{n=-\infty}^0$(resp.
$\{x_n\}_{n=0}^{+\infty}$).
We claim that under the plaque expensiveness assumption the
above sequence
$\{y_i^*\}_{i=-\infty}^{+\infty}$  must be unique. Otherwise, assume
$\{\bar{y}_i^*\}_{i=-\infty}^{+\infty}$ is another $\varepsilon$-quasi-shadowing sequence as in Theorem A in which
$\bar{y}_0^*$ is uniquely determined by $\bar{y}_0^s$ and
$\bar{y}_0^u$ (which are the analogous of  $y_0^s$ and $y_0^u$ to
$y_0^*$). Clearly, $d(y_i^*, \bar{y}_i^*)\leq d(y_i^*, x_i)+d(x_i,
\bar{y}_i^*)<2\varepsilon=\eta$ for any $i\in \mathbf{Z}$. Then by
the plaque expensiveness, each pair of $y_i^*$ and $\bar{y}_i^*$
must lie in a common center plaque. In particular, $y_0^*\in
W_\eta^{c}(\bar{y}_0^*)$. However, since
$\{y_i^*\}_{i=-\infty}^{+\infty}$ is different from
$\{\bar{y}_i^*\}_{i=-\infty}^{+\infty}$ then $\{y_0^s, y_0^u\}$ must
be different from $\{\bar{y}_0^s, \bar{y}_0^u\}$. Remember that
$y_0^s, \bar{y}_0^s\in W_\varepsilon^{s}(x_0)$, $y_0^u,
\bar{y}_0^u\in W_\varepsilon^{u}(x_0)$, $y_0^*\in
W_\varepsilon^{s}(y_0^u)\cap W_\varepsilon^{cu}(y_0^s)$ and
$\bar{y}_0^*\in W_\varepsilon^{s}(\bar{y}_0^u)\cap
W_\varepsilon^{cu}(\bar{y}_0^s)$. Therefore, by the dynamical
coherence, $y_0^*$ and $\bar{y}_0^*$ can not lie in a common center
plaque. A contradiction. Hence the above claim holds. Furthermore, we can see by this claim that $y_0^u$ (resp. $y_0^s$) is the unique
 limit point of the set $\{y_{0,n}^u\}_{n=0}^{\infty}$  (resp. $\{y_{0,n}^s\}_{n=-\infty }^0$).

Let $g$ is $C^0$-close to $f$ with $d(f,g)<\delta$, then for any $x\in
M$, the $g$-orbit of $\{g^i(x)\}_{i=-\infty}^{+\infty}$ is a
$\delta$-pseudo orbit of $f$. Hence there
is a unique sequence $\{y_i^*\}_{i=-\infty}^{+\infty}$ as in Theorem A which
$\varepsilon$-quasi-shadows $\{g^i(x)\}_{i=-\infty}^{+\infty}$.  Let $\pi(x)=y_0^*$,
then by the above claim $\pi$ is a well-defined map on $M$.  Now we show that $\pi$ is continuous. For any $x\in M$ and
any sequence $\{x(k)\}_{k=1}^\infty$ tends to $x$ as $k$ tends to infinity, we get a
sequence of $g$-orbits $\big\{\{g^i(x(k))\}_{i=-\infty}^{+\infty}\big\}_{k=1}^\infty$. Correspondingly, we get a sequence of
quasi-shadowing sequences $\big\{\{y^*_i(k)\}_{i=-\infty}^{+\infty}\big\}_{k=1}^\infty$ in which
$\{y^*_i(k)\}_{i=-\infty}^{+\infty}$ $\varepsilon$-quasi-shadows
$\{g^i(x(k))\}_{i=-\infty}^{+\infty}$. By the definition of
$\pi$, $\pi(x(k))=y^*_0(k)$. Note that for any
$n\in\mathbf{Z}^+$, the $g$-orbit pieces $\{g^i(x(k))\}_{i=-n}^{n}$
converges to $\{g^i(x)\}_{i=-n}^{n}$ (with respect to the Hausdorff
distance) as $k$ tends to infinity. As in Theorem A, we get two sequences
$\{y_{0,n}^u(k)\in W^u_{\varepsilon}(x(k))\}_{k=1}^{\infty}$ and $\{y_{0,-n}^s(k)\in W^s_{\varepsilon}(x(k))\}_{k=1}^{\infty}$. From the construction of
them and the dynamical coherence, it is easy
to see that $y_{0,n}^u(k)\rightarrow y_{0,n}^u$ and $y_{0,-n}^s(k)\rightarrow y_{0,-n}^s$ as
$k\rightarrow +\infty$. From the above claim, for any $k\in \mathbf{Z}^+$ there is also a unique
 limit point $y_0^u(k)$ (resp. $y_0^s(k)$) of the set $\{y_{0,n}^u(k)\}_{n=1}^{\infty}$  (resp. $\{y_{0,-n}^s(k)\}_{n=1}^{\infty}$).  Then we have that
$y_0^{u}(k)\rightarrow y_{0}^{u},y_0^s(k)\rightarrow y_{0}^s$, and hence $\pi(x(k))=y_0^*(k)\rightarrow y_{0}^*=\pi(x)$ as $k\rightarrow +\infty$. This means that $\pi$ is continuous. Obviously, $d(\pi, \text{id}_M)<\varepsilon$,  so from Lemma 3 of \cite{Walters}, $\pi$ is surjective when $\varepsilon$ is small enough.

Now for any $x\in M$, let $\{y_i^*\}_{i=-\infty}^{+\infty}$  be the unique $\varepsilon$-quasi-shadowing sequence of the $g$-orbit $\{g^i(x)\}_{i=-\infty}^{+\infty}$. As in the proof of Theorem A, there is another sequence $\{(y_i^*)'\}_{i=-\infty}^{+\infty}$ is constructed simultaneously. From the construction of these two sequences, we have, for any $i\in \mathbf{Z}$,
$(y_i^*)'=f(y_{i-1}^*)$, and $y_i^*$ lies in the center plaque of $(y_i^*)'$. By uniqueness of these two sequences and the definition of $\pi$, we have $\pi(g^i(x))=y_i^*$. So if we denote $\tau_{g^i(x)}((y_i^*)')=y_i^*$, then we have $$
\pi\circ g^i(x)=\tau_{g^{i-1}(x)}\circ f\circ\pi(g^{i-1}(x)).
$$
This gives (\ref{conjugacy2}). Similar to that for $\pi$, we can get the  continuous dependence of $\{\tau_x: x\in M\}$ on $x$ from their construction.

This completes the proof of Theorem B.

\end{proof}

\end{document}